\documentclass{gtart}

\def\ifplaintex{\expandafter\ifx\csname documentclass\endcsname\relax}

\ifplaintex 
\else
\fi

\expandafter\ifx\csname beginpicture\endcsname\relax
\expandafter\ifx\csname documentclass\endcsname\relax
\input pictex \else
\input prepictex \input pictex \input postpictex \fi\fi

\def\gt{{\mathsurround=0pt\it $\cal G\mskip-2mu$eometry \&\ 
$\cal T\!\!$opology}}        %  journal title in recommended style

\def\gtp{{\mathsurround=0pt\it $\cal G\mskip-2mu$eometry \&\ 
$\cal T\!\!$opology $\cal P\!$ublications}}  % GT publications

\def\lognumber#1{\def\thelognumber{#1}}
\def\volumenumber#1{\def\thevolumenumber{#1}}
\def\papernumber#1{\def\thepapernumber{#1}}
\def\volumeyear#1{\def\thevolumeyear{#1}}

\def\pagenumbers#1#2{\def\startpage{#1}\def\finishpage{#2}}
\def\published#1{\def\publishdate{#1}}
\def\proposed#1{\def\theproposer{#1}}
\def\seconded#1{\def\theseconders{#1}}
\def\received#1{\def\receiveddate{#1}}

\def\accepted#1{\def\accepteddate{#1}}

\long\def\asciiabstract#1{\long\def\theasciiabstract{#1}}

\let\\\par\let\thelognumber\relax
\let\thevolumenumber\relax\let\thepapernumber\relax
\let\thevolumeyear\relax\let\thesamplenumber\relax\let\startpage\relax
\let\finishpage\relax\let\publishdate\relax\let\receiveddate\relax
\let\reviseddate\relax\let\accepteddate\relax\let\theasciititle\relax
\let\theasciiauthors\relax
\let\theasciiabstract\relax
\let\theasciiemail\relax\let\theshortauthors\relax\let\theshorttitle\relax

\long\def\maketitlep{   % start of definition of \maketitlep

\count0=\startpage

\gt\hfill      %   Journal title (top left) 
\beginpicture
\setcoordinatesystem units <0.33truein, 0.33truein> point at 2.2 0.9
\setplotsymbol ({$\cal G$})
\plotsymbolspacing=9truept
\circulararc 315 degrees from 0 1 center at 0 0
\setplotsymbol ({$\cal T$})
\circulararc 315 degrees from 1 -1 center at 1 0
\endpicture
\break
{\small\ifx\thesamplenumber\relax % sample?  
Volume \else Sample
\fi\thevolumenumber\ (\thevolumeyear)
\startpage--\finishpage\nl
Published: \publishdate}
\vglue 0.5truein plus 0.4fil minus 0.1truein

{\parskip=0pt\leftskip 0pt plus 1fil\def\\{\par\smallskip}{\ifplaintex\large
\else\Large\fi\bf\thetitle}\par\medskip}   

\vglue 0pt plus 0.1fil 

{\parskip=0pt\leftskip 0pt plus 1fil\def\\{\par}{\sc\theauthors}
\par\medskip}

\vglue 0pt plus 0.1fil 

{\small\parskip=0pt\let\newline\\
{\leftskip 0pt plus 1fil\def\\{\par}{\sl\theaddress}\par}
\expandafter\ifx\theemail\relax    % email address?
\relax\else\vglue 5pt plus 0.02fil minus 2pt\def\\{\stdspace{\rm 
and}\stdspace} 
\cl{Email:\stdspace\tt\theemail}\fi
\ifx\theurl\relax                  % URL given?
\relax\else\vglue 5pt plus 0.02fil minus 2pt\def\\{\stdspace{\rm 
and}\stdspace}
\cl{URL:\stdspace\tt\theurl}\fi\par}

\vglue 7pt plus 0.3fil minus 3pt

{\bf Abstract}
\vglue 5pt plus 0.1fil minus 2pt

\theabstract

\vglue 7pt plus 0.3fil minus 3pt

{\bf AMS Classification numbers}\quad Primary:\quad \theprimaryclass

Secondary:\quad \thesecondaryclass

\vglue 5pt plus 0.3fil minus 2pt

{\bf Keywords}\quad \thekeywords

\vglue 10pt plus 0.5fil minus 5pt

{\small  Proposed: \theproposer\hfill Received: \receiveddate\nl
Seconded: \theseconders\hfill 
\ifx\reviseddate\relax                         % paper revised?
Accepted: \accepteddate                        % no
\else
Revised: \reviseddate                          % yes
\fi}
\eject
}       %  end of definition of \maketitlep

\let\maketitlepage\maketitlep
\let\maketitle\maketitlepage

\font\phead=cmsl9 scaled 950
\font=cmsl9 scaled 1050
\font\pnum=cmbx10 scaled 913
\font\lnum=cmbx10 
\font\pfoot=cmsl9 scaled 950
\font=cmsl9 scaled 1050
\ifplaintex
\headline{\vbox to 0pt{\vskip -4.5mm\line{\small\phead\ifnum
\count0=\startpage ISSN 1364-0380 (on line)
1465-3060 (printed) \hfill {\pnum\folio}\else\ifodd\count0\def\\{ }% 
\ifx\theshorttitle\relax\thetitle\else\theshorttitle\fi\hfill{\pnum\folio}
\else\def\\{ and }{\pnum\folio}\hfill\ifx\theshortauthors\relax\theauthors
\else\theshortauthors\fi\fi\fi}\vss}}
\footline{\vbox to 0pt{\vglue 0mm\line{\small\pfoot\ifnum\count0=\startpage
\copyright\ \gtp\hfill\else
\gt, Volume \thevolumenumber\ (\thevolumeyear)\hfill\fi}\vss
}}
\else
\makeatletter
\def\@oddhead{{\small\lhead\ifnum\count0=\startpage ISSN 1364-0380 (on line)
1465-3060 (printed) \hfill {\lnum\number\count0}\else\ifodd\count0
\def\\{ }\ifx\theshorttitle\relax \thetitle \else\theshorttitle\fi\hfill
{\lnum\number\count0}\else\def\\{ and }{\lnum\number\count0}
\hfill\ifx\theshortauthors\relax 
\theauthors\else\theshortauthors\fi\fi\fi}}\def\@evenhead{\@oddhead}
\def\@oddfoot{\small\lfoot\ifnum\count0=\startpage\copyright\ \gtp\hfill\else
\gt, Volume \thevolumenumber\ (\thevolumeyear)\hfill\fi}
\def\@evenfoot{\@oddfoot}
\makeatother
\fi

\newwrite\gtoutfile
\long\gdef\makeheadfile{  %%% start of definition of \makeheadfile
{\def\\{, }\def\s{ }
\immediate\openout\gtoutfile head.xxx
\immediate\write\gtoutfile{To: math@arxiv.org}
\immediate\write\gtoutfile{Subject: put or rep NNNNN:pppp}
\immediate\write\gtoutfile{--text follows this line--}
\immediate\write\gtoutfile{Proxy-for: \ifx\theasciiauthors\relax
\theauthors\else\theasciiauthors\fi\s<\ifx\theasciiemail\relax\theemail\else\theasciiemail\fi>}
\immediate\write\gtoutfile{\noexpand\\}
\immediate\write\gtoutfile{Authors: \ifx\theasciiauthors\relax
\theauthors\else\theasciiauthors\fi}
{\def\\{ }\immediate\write\gtoutfile{Title: \ifx\theasciititle\relax
\thetitle\else\theasciititle\fi}}
\immediate\write\gtoutfile{Subj-class: GT or SG or MG etc}
\immediate\write\gtoutfile{MSC-class: \theprimaryclass\ifx\thesecondaryclass\relax\else, \thesecondaryclass\fi}
\immediate\write\gtoutfile{Journal-ref: Geom. Topol. \thevolumenumber
(\thevolumeyear) \startpage-\finishpage}
\immediate\write\gtoutfile{Comments: Published by Geometry and Topology at}
\immediate\write\gtoutfile{\s\s http://www.maths.warwick.ac.uk/gt/GTVol\thevolumenumber/paper\thepapernumber.abs.html}
\immediate\write\gtoutfile{\noexpand\\}
\immediate\write\gtoutfile{}
\ifx\theasciiabstract\relax
\immediate\write\gtoutfile{\theabstract}\else
\immediate\write\gtoutfile{\theasciiabstract}\fi
\immediate\write\gtoutfile{}
\immediate\write\gtoutfile{\noexpand\\}
\immediate\write\gtoutfile{}
\immediate\closeout\gtoutfile}}  %%% end of definition of \makeheadfile

\def\maketitlepage{\maketitlep\makeheadfile}
\let\maketitle\maketitlepage

\def\ifplaintex{\expandafter\ifx\csname documentclass\endcsname\relax}

\ifplaintex 
\else
\fi

\expandafter\ifx\csname beginpicture\endcsname\relax
\expandafter\ifx\csname documentclass\endcsname\relax
\input pictex \else
\input prepictex \input pictex \input postpictex \fi\fi

\def\gt{{\mathsurround=0pt\it $\cal G\mskip-2mu$eometry \&\ 
$\cal T\!\!$opology}}        %  journal title in recommended style

\def\gtp{{\mathsurround=0pt\it $\cal G\mskip-2mu$eometry \&\ 
$\cal T\!\!$opology $\cal P\!$ublications}}  % GT publications

\def\lognumber#1{\def\thelognumber{#1}}
\def\volumenumber#1{\def\thevolumenumber{#1}}
\def\papernumber#1{\def\thepapernumber{#1}}
\def\volumeyear#1{\def\thevolumeyear{#1}}

\def\pagenumbers#1#2{\def\startpage{#1}\def\finishpage{#2}}
\def\published#1{\def\publishdate{#1}}
\def\proposed#1{\def\theproposer{#1}}
\def\seconded#1{\def\theseconders{#1}}
\def\received#1{\def\receiveddate{#1}}

\def\accepted#1{\def\accepteddate{#1}}

\long\def\asciiabstract#1{\long\def\theasciiabstract{#1}}

\let\\\par\let\thelognumber\relax
\let\thevolumenumber\relax\let\thepapernumber\relax
\let\thevolumeyear\relax\let\thesamplenumber\relax\let\startpage\relax
\let\finishpage\relax\let\publishdate\relax\let\receiveddate\relax
\let\reviseddate\relax\let\accepteddate\relax\let\theasciititle\relax
\let\theasciiauthors\relax
\let\theasciiabstract\relax
\let\theasciiemail\relax\let\theshortauthors\relax\let\theshorttitle\relax

\long\def\maketitlep{   % start of definition of \maketitlep

\count0=\startpage

\gt\hfill      %   Journal title (top left) 
\beginpicture
\setcoordinatesystem units <0.33truein, 0.33truein> point at 2.2 0.9
\setplotsymbol ({$\cal G$})
\plotsymbolspacing=9truept
\circulararc 315 degrees from 0 1 center at 0 0
\setplotsymbol ({$\cal T$})
\circulararc 315 degrees from 1 -1 center at 1 0
\endpicture
\break
{\small\ifx\thesamplenumber\relax % sample?  
Volume \else Sample
\fi\thevolumenumber\ (\thevolumeyear)
\startpage--\finishpage\nl
Published: \publishdate}
\vglue 0.5truein plus 0.4fil minus 0.1truein

{\parskip=0pt\leftskip 0pt plus 1fil\def\\{\par\smallskip}{\ifplaintex\large
\else\Large\fi\bf\thetitle}\par\medskip}   

\vglue 0pt plus 0.1fil 

{\parskip=0pt\leftskip 0pt plus 1fil\def\\{\par}{\sc\theauthors}
\par\medskip}

\vglue 0pt plus 0.1fil 

{\small\parskip=0pt\let\newline\\
{\leftskip 0pt plus 1fil\def\\{\par}{\sl\theaddress}\par}
\expandafter\ifx\theemail\relax    % email address?
\relax\else\vglue 5pt plus 0.02fil minus 2pt\def\\{\stdspace{\rm 
and}\stdspace} 
\cl{Email:\stdspace\tt\theemail}\fi
\ifx\theurl\relax                  % URL given?
\relax\else\vglue 5pt plus 0.02fil minus 2pt\def\\{\stdspace{\rm 
and}\stdspace}
\cl{URL:\stdspace\tt\theurl}\fi\par}

\vglue 7pt plus 0.3fil minus 3pt

{\bf Abstract}
\vglue 5pt plus 0.1fil minus 2pt

\theabstract

\vglue 7pt plus 0.3fil minus 3pt

{\bf AMS Classification numbers}\quad Primary:\quad \theprimaryclass

Secondary:\quad \thesecondaryclass

\vglue 5pt plus 0.3fil minus 2pt

{\bf Keywords}\quad \thekeywords

\vglue 10pt plus 0.5fil minus 5pt

{\small  Proposed: \theproposer\hfill Received: \receiveddate\nl
Seconded: \theseconders\hfill 
\ifx\reviseddate\relax                         % paper revised?
Accepted: \accepteddate                        % no
\else
Revised: \reviseddate                          % yes
\fi}
\eject
}       %  end of definition of \maketitlep

\let\maketitlepage\maketitlep
\let\maketitle\maketitlepage

\font\phead=cmsl9 scaled 950
\font=cmsl9 scaled 1050
\font\pnum=cmbx10 scaled 913
\font\lnum=cmbx10 
\font\pfoot=cmsl9 scaled 950
\font=cmsl9 scaled 1050
\ifplaintex
\headline{\vbox to 0pt{\vskip -4.5mm\line{\small\phead\ifnum
\count0=\startpage ISSN 1364-0380 (on line)
1465-3060 (printed) \hfill {\pnum\folio}\else\ifodd\count0\def\\{ }% 
\ifx\theshorttitle\relax\thetitle\else\theshorttitle\fi\hfill{\pnum\folio}
\else\def\\{ and }{\pnum\folio}\hfill\ifx\theshortauthors\relax\theauthors
\else\theshortauthors\fi\fi\fi}\vss}}
\footline{\vbox to 0pt{\vglue 0mm\line{\small\pfoot\ifnum\count0=\startpage
\copyright\ \gtp\hfill\else
\gt, Volume \thevolumenumber\ (\thevolumeyear)\hfill\fi}\vss
}}
\else
\makeatletter
\def\@oddhead{{\small\lhead\ifnum\count0=\startpage ISSN 1364-0380 (on line)
1465-3060 (printed) \hfill {\lnum\number\count0}\else\ifodd\count0
\def\\{ }\ifx\theshorttitle\relax \thetitle \else\theshorttitle\fi\hfill
{\lnum\number\count0}\else\def\\{ and }{\lnum\number\count0}
\hfill\ifx\theshortauthors\relax 
\theauthors\else\theshortauthors\fi\fi\fi}}\def\@evenhead{\@oddhead}
\def\@oddfoot{\small\lfoot\ifnum\count0=\startpage\copyright\ \gtp\hfill\else
\gt, Volume \thevolumenumber\ (\thevolumeyear)\hfill\fi}
\def\@evenfoot{\@oddfoot}
\makeatother
\fi

\newwrite\gtoutfile
\long\gdef\makeheadfile{  %%% start of definition of \makeheadfile
{\def\\{, }\def\s{ }
\immediate\openout\gtoutfile head.xxx
\immediate\write\gtoutfile{To: math@arxiv.org}
\immediate\write\gtoutfile{Subject: put or rep NNNNN:pppp}
\immediate\write\gtoutfile{--text follows this line--}
\immediate\write\gtoutfile{Proxy-for: \ifx\theasciiauthors\relax
\theauthors\else\theasciiauthors\fi\s<\ifx\theasciiemail\relax\theemail\else\theasciiemail\fi>}
\immediate\write\gtoutfile{\noexpand\\}
\immediate\write\gtoutfile{Authors: \ifx\theasciiauthors\relax
\theauthors\else\theasciiauthors\fi}
{\def\\{ }\immediate\write\gtoutfile{Title: \ifx\theasciititle\relax
\thetitle\else\theasciititle\fi}}
\immediate\write\gtoutfile{Subj-class: GT or SG or MG etc}
\immediate\write\gtoutfile{MSC-class: \theprimaryclass\ifx\thesecondaryclass\relax\else, \thesecondaryclass\fi}
\immediate\write\gtoutfile{Journal-ref: Geom. Topol. \thevolumenumber
(\thevolumeyear) \startpage-\finishpage}
\immediate\write\gtoutfile{Comments: Published by Geometry and Topology at}
\immediate\write\gtoutfile{\s\s http://www.maths.warwick.ac.uk/gt/GTVol\thevolumenumber/paper\thepapernumber.abs.html}
\immediate\write\gtoutfile{\noexpand\\}
\immediate\write\gtoutfile{}
\ifx\theasciiabstract\relax
\immediate\write\gtoutfile{\theabstract}\else
\immediate\write\gtoutfile{\theasciiabstract}\fi
\immediate\write\gtoutfile{}
\immediate\write\gtoutfile{\noexpand\\}
\immediate\write\gtoutfile{}
\immediate\closeout\gtoutfile}}  %%% end of definition of \makeheadfile

\def\maketitlepage{\maketitlep\makeheadfile}
\let\maketitle\maketitlepage

\lognumber{210}

\volumenumber{6}\papernumber{10}\volumeyear{2002}
\pagenumbers{329}{353}
\received{5 October 2001}
\accepted{25 May 2002}
\published{25 June 2002}

\proposed{Robion Kirby}
\seconded{Michael Freedman, Walter Neumann}

\usepackage{amsmath,amssymb}
\def\S{Section }

\newcommand{\R}{\mathbb{R}}
\newcommand{\Z}{\mathbb{Z}}
\newcommand{\C}{\mathbb{C}}
\newcommand{\Imag}{\mathop{\textup{Im}}}
\newcommand{\Vol}{\mathop{\textup{vol}}}
\newcommand{\supp}{\mathop{\textup{supp}}}
\newcommand{\newchi}{\chi^{\phantom0}}
\newcommand{\newchihat}{\widehat{\chi}^{\phantom0}}

\theoremstyle{plain}
\newtheorem{theorem}{Theorem}[section]
\newtheorem{corollary}[theorem]{Corollary}
\newtheorem{lemma}[theorem]{Lemma}
\newtheorem{proposition}[theorem]{Proposition}

\theoremstyle{definition}
\newtheorem{definition}[theorem]{Definition}
\newtheorem{conjecture}[theorem]{Conjecture}
\newtheorem{question}[theorem]{Open Question}

\numberwithin{equation}{section}

\title{New upper bounds on sphere packings II}
\author{Henry Cohn}
\address{Microsoft Research, One Microsoft Way\\Redmond, WA 98052-6399, USA}
\email{cohn@microsoft.com}
\url{http://research.microsoft.com/\char'176 cohn}


\begin{document}

\begin{abstract}
We continue the study of the linear programming bounds for sphere
packing introduced by Cohn and Elkies.  We use theta series to
give another proof of the principal theorem, and present some
related results and conjectures.

This article is in the arXiv as: {\tt arXiv:math.MG/0110010}
\bigskip
\begin{center}
\textit{Dedicated to Daniel Lewin}\qua (14 May 1970 -- 11 September 2001)
\end{center}
\end{abstract}
\asciiabstract{%
We continue the study of the linear programming bounds for sphere
packing introduced by Cohn and Elkies.  We use theta series to
give another proof of the principal theorem, and present some
related results and conjectures.}
\maketitlepage

\section{Introduction}

In \cite{CE}, Cohn and Elkies introduce linear programming bounds
for the sphere packing problem, and use them to prove new upper
bounds on the sphere packing density in low dimensions.  These
bounds are the best bounds known in dimensions $4$ through $36$,
and seem to be sharp in dimensions $8$ and $24$, although that has
not yet been proved.  Here, we continue the study of these bounds,
by giving another derivation of the main theorem of \cite{CE}.  We
then prove an optimality theorem of Gorbachev \cite{G}, and
outline in some conjectures how the proof techniques should apply
more generally.

We continue to use the notation of \cite{CE}.  See the introduction
of that paper for background and references on sphere packing.

The main theorem Cohn and Elkies prove is the following:

\begin{theorem}
\label{main}
Suppose $f \co \R^n \rightarrow \R$ is a radial,
admissible function, is not identically zero, and satisfies the
following two conditions:
\begin{enumerate}
\item $f(x)\le 0$ for $|x| \ge 1$, and

\item ${\widehat f}(t)\ge 0$ for all $t$.
\end{enumerate}
Then the center densities of $n$--dimensional sphere packings are
bounded above by
$$
\frac{f(0)}{2^n{\widehat f}(0)}.
$$
\end{theorem}

Here, the Fourier transform is normalized by
$$
\widehat{f}(t) =
\int_{\R^n} f(x) e^{2\pi i \langle t,x \rangle} \, dx,
$$
and \textit{admissibility\/} means that there is a constant
$\varepsilon>0$ such that both $|f(x)|$ and $|\widehat{f}(x)|$
are bounded above by a constant times $(1+|x|)^{-n-\varepsilon}$.
More broadly, we could in fact take $f$ to be any function to
which the Poisson summation formula applies: for every lattice
$\Lambda \subset \R^n$ and every vector $v \in \R^n$,
$$
\sum_{x \in \Lambda} f(x+v) = \frac{1}{\Vol(\R^n/\Lambda)}
\sum_{t \in \Lambda^*} e^{-2\pi i \langle v,t \rangle}
\widehat{f}(t).
$$
However, the narrower definition of admissibility is easier to
check and seemingly suffices for all natural examples.

Section~\ref{newproof} gives another proof of Theorem~\ref{main},
for $n>1$.  This proof is not as simple as the one in \cite{CE},
but the method is of interest in its own right, as are some of the
intermediate results. Section~\ref{opt} proves Gorbachev's
theorem \cite{G} that certain admissible functions (those
constructed in Proposition~6.1
of \cite{CE}, or independently by Gorbachev) are optimal, among
functions whose Fourier transforms have support in a certain
ball.  Finally, Section~\ref{dualprog} discusses the dual linear
program, and puts the techniques of Section~\ref{opt} into a
broader context.

\rk{Acknowledgements} 

I thank Richard Askey, Noam Elkies, Pavel Etingof, and George Gasper
for their advice about special functions, and Tom Brennan, Harold
Diamond, Gerald Folland, Cormac Herley, David Jerison, Greg Kuperberg,
Ben Logan, L\'aszl\'o Lov\'asz, Steve Miller, Amin Shokrollahi, Neil
Sloane, Hart Smith, Jeffrey Vaaler, David Vogan, and Michael Weinstein
for helpful discussions. I was previously supported by an NSF Graduate
Research Fellowship and a summer internship at Lucent Technologies,
and currently hold an American Institute of Mathematics Five-year
Fellowship.

\section{Positivity of theta series coefficients}
\label{newproof}

We will prove Theorem~\ref{main} using the positivity of the
coefficients of the theta series of lattices. For each lattice,
the theta series of its dual must have positive coefficients, and
these coefficients are some transformation of those for the
original lattice.  This puts strong constraints on the theta
series of a lattice, which we exploit below.  For simplicity, we
will deal only with the case of lattice packings, but everything
in this section applies to all sphere packings, by replacing the
theta series of a lattice with the average theta series of a
periodic packing (see \cite[page 45]{SPLAG}). Also, for technical
reasons we will deal only with the case $n>1$, which is not a
serious restriction as $1$--dimensional sphere packing is trivial.

Unfortunately, carrying this program out rigorously involves
dealing with a number of technicalities.  If one simply wants an
idea of the overall argument, without worrying about rigor, one
can follow this plan:  Ignore Lemma~\ref{fejer} and all
references to Ces\`aro sums, and assume that all Laguerre series
converge. Ignore the uniformity of convergence in
Lemma~\ref{saddlept} (in which case the proof becomes far
simpler). Ignore the justification of interchanging the sum and
integral in Lemma~\ref{sumint}.  Following this plan will of
course not lead to a rigorous proof, but it may make the
underlying ideas clearer.

Before going further, we need a lemma about Laguerre polynomials.
Let $L^\alpha_k$ be the Laguerre polynomial of degree $k$ and
parameter $\alpha > -1$. These polynomials are orthogonal with
respect to the weight $x^{\alpha} e^{-x}\, dx$ on $[0,\infty)$.

\begin{lemma}
\label{induction}
For every non-negative integer $k$, $\alpha >
-1$, and $y \in \R$, we have
$$
\frac{(-1)^k}{k!} \frac{d^k}{du^k} \left( u^{-\alpha-1} e^{-y/u}
\right) = u^{-\alpha-1-k} e^{-y/u} L_k^\alpha(y/u).
$$
\end{lemma}

\begin{proof}
This is easily proved by induction, using standard  properties of
Laguerre polynomials (see Section~6.2 of \cite{AAR}, or
Sections~4.17--4.24 of \cite{Leb}).
\end{proof}

Suppose $\Lambda \subset \R^n$ is a lattice, and define a measure
$\mu$ on $[0,\infty)$ consisting of a point mass at $x$ for each
vector in $\Lambda$ of norm $x$, where the norm of $v$ is
$\langle v, v \rangle$. The purpose of $\mu$ is to allow us to
sum over all lattice vectors without having to index the sum in
our notation;  instead, we simply integrate with respect to $\mu$.
Although $\mu$ depends on $\Lambda$, for simplicity our notation
does not make that dependence explicit.

The key positivity property of $\mu$ is the following lemma:

\begin{lemma}
\label{thetacoeff}
For all $y > 0$ and all non-negative integers $k$,
$$
\int_0^\infty L^{n/2-1}_k(xy) e^{-xy} \, d\mu(x) \ge 0.
$$
\end{lemma}

\begin{proof}
The theta series of $\Lambda$ is given by
$$
\Theta_\Lambda(z) = \int_0^\infty e^{i\pi x z} \,d\mu(x),
$$
and it follows from the Poisson summation formula that the theta
series of the dual lattice $\Lambda^*$ is given by
$$
\Theta_{\Lambda^*}(z) = \Vol(\R^n/\Lambda)
\left(\frac{i}{z}\right)^{n/2} \Theta_\Lambda\left(
-\frac{1}{z}\right).
$$
(See equation~(19) in \cite[page 103]{SPLAG}.)

It will be more convenient for us to work with the variable $y$
given by $y = - i \pi z$.  Let $T(y) = \Theta_\Lambda(z)$, so that
$$
T(y) = \int_0^\infty e^{-xy} \,d\mu(x).
$$
Then up to a positive factor, the theta series of $\Lambda^*$
is given by $y^{-n/2} T(\pi^2/y)$.

We know that $y^{-n/2} T(\pi^2/y)$ is a positive linear
combination of functions $e^{-cy}$ with $c \ge 0$, because it is
the theta series of a lattice (times a positive constant). Hence,
its successive derivatives with respect to $y$ alternate in sign.
We have
$$
y^{-n/2} T(\pi^2/y) = \int_0^\infty y^{-n/2} e^{-\pi^2 x/y} \,
d\mu(x),
$$
{}from which it follows using Lemma~\ref{induction} that
$$
\frac{(-1)^k}{k!} \frac{d^k}{dy^k} \left( y^{-n/2} T(\pi^2/y)
\right) = \int_0^\infty y^{-n/2-k} e^{-\pi^2 x/y}
L_k^{n/2-1}(\pi^2 x/y)\, d\mu(x).
$$
(Differentiating under the integral sign,
which really denotes a sum,
is justified by uniform convergence of the
differentiated sum; see Theorem~7.17 of \cite{Ru}.) Now the
change of variable $y \leftrightarrow \pi^2/y$ shows us that
$$
\int_0^\infty L^{n/2-1}_k(xy) e^{-xy} \, d\mu(x) \ge 0,
$$
as desired.
\end{proof}

When we use only the fact that the derivatives of
$y^{-n/2}T(\pi^2/y)$ alternate in sign, we do not lose much
information---by a theorem of Bernstein (see \S12 of Chapter~IV
of \cite{W}), this property characterizes functions of the form
$$
\int_0^\infty e^{-xy} \, d\mu(x)
$$
for some measure $\mu$ on
$[0,\infty)$. Also, it is not surprising that the inequalities in
Lemma~\ref{thetacoeff} occur for all scalings $y$, because so far
our setup is scale-invariant.

If the shortest non-zero vectors in $\Lambda$ have length $1$
(that is, $\Lambda$ leads to a packing with balls of radius $1/2$),
then the center density of the lattice packing given by $\Lambda$
equals
$$
(4\pi)^{-n/2}\lim_{y \rightarrow 0+} y^{n/2}T(y).
$$
The proof is as follows.  The relationship between the theta
series of $\Lambda^*$ and $\Lambda$ is
$$
\frac{T_{\Lambda^*}(y)}{2^n\Vol(\R^n/\Lambda)} = (4\pi)^{-n/2}
\left(\frac{\pi^2}{y}\right)^{n/2}
T_\Lambda\left(\frac{\pi^2}{y}\right).
$$
As we let $y \rightarrow \infty$, the right hand side becomes the
limit above, and the left hand side tends to $1/(2^n
\Vol(\R^n/\Lambda))$, which is the center density.

Using Lemma~\ref{thetacoeff}, we can bound the center density.
First, we need a definition and a lemma.

\begin{definition}
\label{silpdef}
A function $f \co [0,\infty) \to \R$ has the $\alpha$--SILP property
(``scale-invariant Laguerre positivity'') if the following
conditions hold:
\begin{enumerate}
\item $f$ is continuous and for some $\varepsilon > 0$ and $C>0$,
we have
$$
|f(x)| \le C(1+|x|)^{-\alpha-1-\varepsilon}
$$
for all $x$, and
\item for every $y>0$, the Laguerre series
$$
\sum_{j \ge 0} a_j(y) L_j^\alpha(x),
$$
for $x \mapsto f(x/y)$ has $a_j(y) \ge 0$ for all $j$.
\end{enumerate}
\end{definition}

Condition~(1) is merely a technical restriction; condition~(2) is
the heart of the matter.  Notice that the orthogonality of the
Laguerre polynomials implies that
$$
a_j(y) = \frac{\int_0^\infty f(x/y)L^{\alpha}_j(x) x^{\alpha}
e^{-x} \, dx} {\int_0^\infty L^{\alpha}_j(x)^2 x^{\alpha} e^{-x}
\, dx} = \frac{\int_0^\infty f(x/y)L^{\alpha}_j(x) x^{\alpha}
e^{-x} \, dx} {\Gamma(j+\alpha+1)/j!}.
$$
We make no assumption about convergence for the Laguerre series
in Definition~\ref{silpdef}.
However, the following
analogue of Fej\'er's theorem on Fourier series
holds.  It is a
simple consequence of results in \cite{T}.  We could also make
use of \cite{Stem} to prove a marginally weaker result (which
would still suffice for our purposes).

\begin{lemma}
\label{fejer}
Let $\alpha \ge 0$, and let $f \co [0,\infty) \to \R$
be an $\alpha$--SILP function.  Then for all $k > \alpha + 1/2$,
the $(C,k)$ Ces\`aro means
$$
\binom{k+m}{m}^{-1} \sum_{j=0}^m \binom{k+m-j}{m-j} a_j(y)
L_j^\alpha(x) e^{-x/2}
$$
of the partial sums of the  series
$$
\sum_{j \ge 0} a_j(y) L_j^\alpha(x) e^{-x/2}
$$
converge uniformly to $f(x/y)e^{-x/2}$ on $[0,\infty)$,
as $m \to \infty$. (Here,
$a_j(y)$ is as above.)
\end{lemma}

\begin{proof}
We take $y=1$ for notational simplicity; of course,
the same proof holds for
each $y>0$. For a function $g \co [0,\infty) \to \R$, let $\tilde
g(x) = g(x) e^{-x/2}$, and let $\sigma_m g(x)$ denote the
Ces\`aro mean
$$
\sigma_m g(x) = \binom{k+m}{m}^{-1} \sum_{j=0}^m
\binom{k+m-j}{m-j} b_j L_j^\alpha(x) e^{-x/2},
$$
where $g$ has Laguerre coefficients $b_j$.
Theorem~6.2.1 of \cite{T} says that there exists a constant
$C$ such that for all $m$ and all $g$ such that $\tilde g \in
L^\infty([0,\infty),x^\alpha \, dx)$,
$$
|| \sigma_m g ||_\infty \le C || \tilde g ||_\infty,
$$
where $||\cdot||_\infty$ denotes the norm on
$L^\infty([0,\infty),x^\alpha \, dx)$.

We can then imitate the proof of Theorem~2 in \cite{P}. Let
$\varepsilon>0$. By Theorem~18 of \cite{St}, $\tilde f$ can be
uniformly approximated on $[0,\infty)$ by $\tilde g$ with $g$ a
polynomial. Choose $g$ so that
$$
\left|\left| \tilde f - \tilde g\right|\right|_\infty <
\frac{\varepsilon}{2+2C}.
$$
Then
$$
|| \sigma_m f - \sigma_m g||_\infty < \frac{C\varepsilon}{2+2C}.
$$
For sufficiently large $m$, we have
$$
|| \sigma_m g - \tilde g||_\infty < \frac{\varepsilon}{2},
$$
since $g$ is a polynomial.  It follows that
$$
\left|\left| \sigma_m f - \tilde f \right|\right|_\infty <
\varepsilon.
$$
Thus, $\sigma_m f$ tends uniformly to $f$ as $m \to \infty$.
\end{proof}

Of course, this proof made no use of the positivity of the
Laguerre coefficients, and in fact could be carried out with far
weaker constraints on the behavior of $f$ at infinity.  We stated
it in terms of $\alpha$--SILP functions only because those are the
functions to which we will apply it.  The requirement that
$\alpha$ be non-negative is part of the hypotheses of
Theorem~6.2.1 of \cite{T}.  Perhaps one could prove an analogue
of Lemma~\ref{fejer} for $\alpha < 0$, but in terms of sphere
packing that would cover only the one-dimensional case.

\begin{theorem}
\label{silpbound}
Let $n > 1$.
Suppose $f$ has the $(n/2-1)$--SILP property,
with $f(0)=1$ and $f(x) \le 0$ for $x \ge 1$.  Then the center
density for $n$--dimensional lattice packings is bounded above
by
$$
\frac{\Gamma(n/2)}{2^n \pi^{n/2} \int_0^\infty f(x) x^{n/2-1} \,
dx}.
$$
\end{theorem}

As was pointed out above, the same bound holds for all sphere
packings, not just lattice packings.  One can prove this more
general result by replacing the theta series of a lattice with
the averaged theta series of a periodic packing in
Lemma~\ref{thetacoeff}, but for simplicity we restrict our
attention to lattices.

\begin{proof}
Without loss of generality, we can assume that our lattice is
scaled so as to have packing radius $1/2$ (that is, every non-zero
vector has norm at least $1$).  Define $\mu$, $T$, $a_k(y)$,
etc.\ as before.

We have
$$
f(0) \ge \int_0^\infty f(x) e^{-xy} \, d\mu(x),
$$
since all contributions to the integral from $x > 0$ are
non-positive.

Let $k > (n-1)/2$, and
$$
\sigma_m f(x) = \binom{k+m}{m}^{-1} \sum_{j=0}^m
\binom{k+m-j}{m-j} a_j(y) L_j^{n/2-1}(xy) e^{-xy/2}.
$$
Then
$$
\int_0^\infty \sigma_m f(x) e^{-xy/2} \, d\mu(x) \ge a_0(y)
\int_0^\infty e^{-xy} \, d\mu(x) = a_0(y) T(y),
$$
since by Lemma~\ref{thetacoeff} all the terms in $\sigma_m f(x)$
with $j>0$ contribute a non-negative amount.  Since $\sigma_m
f(x)$ converges uniformly to $f(x) e^{-xy/2}$ as $m \to \infty$ by
Lemma~\ref{fejer} (and because constant functions are integrable
with respect to $e^{-xy/2} \, d\mu(x)$), we have
$$
\lim_{m \to \infty}
\int_0^\infty \sigma_m f(x) e^{-xy/2} \, d\mu(x)
=
\int_0^\infty f(x) e^{-xy} \, d\mu(x).
$$
It follows that
$$
\int_0^\infty f(x) e^{-xy} \, d\mu(x) \ge a_0(y) T(y),
$$
and hence
$$
f(0) \ge a_0(y) T(y).
$$
Thus, the center density
is bounded above by
$$
\lim_{y \rightarrow 0+} \frac{y^{n/2}f(0)}{(4\pi)^{n/2} a_0(y)}.
$$
We can evaluate that limit, since
$$
a_0(y) = \frac{\int_0^\infty f(x/y) x^{n/2-1} e^{-x} \, dx}
{\Gamma(n/2)} = \frac{y^{n/2}\int_0^\infty f(u) u^{n/2-1} e^{-yu}
\, du}{\Gamma(n/2)},
$$
and $\int_0^\infty f(u) u^{n/2-1} e^{-yu} \, du$ converges to
$\int_0^\infty f(u) u^{n/2-1} \, du$ as $y \to 0+$, by dominated
convergence. Applying this formula leads to the bound in the
theorem statement.
\end{proof}

Theorem~\ref{silpbound} amounts to essentially the same bound as
Theorem~\ref{main}, although that is not immediately obvious. The
key is Proposition~\ref{silpchar}, which tells us that there is
essentially only one $\alpha$--SILP function for each $\alpha$, in
the sense that every $\alpha$--SILP function is a positive
combination of scalings of this function.  First, we need two
technical lemmas.

\begin{lemma}
\label{saddlept}
For $\alpha > -1/2$ and $x \in [0,\infty)$,
$$
\lim_{k \to \infty} k^{-\alpha} L_k^\alpha(x/k) e^{-x/k}
= x^{-\alpha/2} J_\alpha(2\sqrt{x}),
$$
and convergence is uniform over $[0,\infty)$.
\end{lemma}

Note that uniform convergence is false for $\alpha=-1/2$, because
$k^{-\alpha} L_k^\alpha(x/k) e^{-x/k}$ tends to~$0$ as $x \to
\infty$ but the right side does not. Since we take $\alpha =
n/2-1$ in dimension~$n$, the only case this rules out is the
trivial $1$--dimensional case, and that is hardly a problem since
it is already ruled out by Theorem~\ref{silpbound} (via
Lemma~\ref{fejer}).

\begin{proof}
Pointwise convergence is known (see 10.12 (36) in
\cite[page 191]{BMP}), but the statements the author knows of
in the literature omit the $e^{-x/k}$ factor that makes the
convergence uniform.

We consider two cases. In the first, $x \ge k^{1+\delta}$ with
$\delta>0$ fixed as $k \to \infty$. Then $x^{-\alpha/2}
J_\alpha(2\sqrt{x})$ tends uniformly to~$0$ as $k \to \infty$,
and we just need to verify that $k^{-\alpha} L_k^\alpha(x/k)
e^{-x/k}$ does as well.  For that, we use Theorem~8.91.2 from
\cite{Sz2}.  It implies that for $a>0$
$$
\max_{x \ge a} \left|e^{-x/2} L_k^\alpha(x)\right| = O(k^C),
$$
where $C = \max(-1/3,\alpha/2-1/4)$.  It follows that
$k^{-\alpha} L_k^\alpha(x/k) e^{-x/k}$ tends uniformly to~$0$ as
$k \to \infty$ with $x \ge k^{1+\delta}$.

Thus, we need only deal with the case of $x \le k^{1+\delta}$.  We
start with (4.19.3) from \cite{Leb} (which holds for all $\alpha >
-1$, not just $\alpha > 1$ as inadvertently stated in \cite{Leb}),
which says that
$$
L_k^\alpha(x) = \frac{e^x x^{-\alpha/2}}{k!}
\int_0^\infty t^{k+\alpha/2} J_\alpha(2\sqrt{xt}) e^{-t} \, dt.
$$
Thus,
\begin{eqnarray*}
k^{-\alpha} L_k^\alpha(x/k) e^{-x/k} &=& \frac{x^{-\alpha/2}
k^{k+1}}{k!} \int_0^\infty t^{\alpha/2} J_\alpha(2\sqrt{xt})
e^{k(\log t - t)} \, dt\\
&=& (1+o(1)) e^k \sqrt{\frac{k}{2\pi}} \int_0^\infty
(t/x)^{\alpha/2} J_\alpha(2\sqrt{xt})
e^{k(\log t - t)} \, dt.
\end{eqnarray*}
The exponent $\log t - t$ is maximized at $t=1$, so we can use the
Laplace method to estimate this integral (see Chapter~4 of
\cite{dB}). In the following calculations, all constants implicit
in big-$O$ terms are independent of $x$.

Let $\varepsilon>0$ be small ($\varepsilon$ will be a function of
$k$). Our integral nearly equals that over the interval
$[1-\varepsilon,1+\varepsilon]$, since for any $C<1/2$ we have
$\log t - t < -1-C\varepsilon^2$ outside
$[1-\varepsilon,1+\varepsilon]$ for sufficiently small
$\varepsilon$, and hence
$$
\left|\int_0^\infty (t/x)^{\alpha/2} J_\alpha(2\sqrt{xt})
e^{k(\log t - t)} \, dt - \int_{1-\varepsilon}^{1+\varepsilon}
(t/x)^{\alpha/2} J_\alpha(2\sqrt{xt}) e^{k(\log t - t)} \, dt
\right|
$$
is bounded by
$$
e^{-(k-1)(1+C \varepsilon^2)} \int_0^\infty t^{\alpha/2}
\left|\frac{J_\alpha(2\sqrt{xt})}{x^{\alpha/2}}\right| e^{\log t
- t} \, dt = O\left(e^{-k(1+C \varepsilon^2)}\right).
$$
Thus, we just need to estimate
$$
\int_{1-\varepsilon}^{1+\varepsilon} (t/x)^{\alpha/2}
J_\alpha(2\sqrt{xt}) e^{k(\log t - t)} \, dt.
$$
We would like to approximate it with
$$
x^{-\alpha/2} J_\alpha(2\sqrt{x})
\int_{1-\varepsilon}^{1+\varepsilon} e^{k(\log t - t)} \, dt.
$$
The difference between these integrals is bounded by a constant
times the product of $\varepsilon$, the maximum of the
$t$--derivative of $(t/x)^{\alpha/2} J_\alpha(2\sqrt{xt})$ over $t
\in [1-\varepsilon,1+\varepsilon]$, and
$$
\int_{1-\varepsilon}^{1+\varepsilon} e^{k(\log t - t)} \, dt.
$$
We have
$$
\frac{\partial}{\partial t} (t^{\alpha/2}J_\alpha(2\sqrt{xt})) =
\frac{\alpha}{2} t^{\alpha/2-1}J_\alpha(2\sqrt{xt}) +
\left(-J_{\alpha+1}(2\sqrt{xt}) + \frac{\alpha
J_\alpha(2\sqrt{xt})}
{2\sqrt{xt}}\right)\frac{t^{\alpha/2}x}{\sqrt{xt}}.
$$
For $x$ near $0$, $x^{-\alpha/2} \partial
(t^{\alpha/2}J_\alpha(2\sqrt{xt}))/\partial t$ remains bounded;
for $x$ away from $0$ it is at most $O(x^{1/4-\alpha/2}),$ which
is at most $O(x^{1/2-\delta})$ if $\delta$ is small
enough relative to $\alpha$ (which we can assume).
Because $x \le k^{1+\delta}$, we have $x^{1/2-\delta}
\le k^{1/2-\delta/2}$.

Thus,
$$
\int_{1-\varepsilon}^{1+\varepsilon} (t/x)^{\alpha/2}
J_\alpha(2\sqrt{xt}) e^{k(\log t - t)} \, dt
$$
equals
$$
\left(x^{-\alpha/2}J_\alpha(2\sqrt{x}) + O\left(\varepsilon
k^{1/2-\delta/2}\right)\right)
\int_{1-\varepsilon}^{1+\varepsilon} e^{k(\log t - t)} \, dt.
$$
If we expand $\log t - t= -1 - (t-1)^2/2 + O((t-1)^3)$, we find
that
$$
\int_{1-\varepsilon}^{1+\varepsilon} e^{k(\log t - t)} \, dt =
(1+ o(1))e^{-k}\sqrt\frac{2\pi}{k},
$$
as long as $k\varepsilon^2 \to \infty$, so that the interval we
are integrating over is much wider than the standard deviation of
the Gaussian we are using to approximate the integrand.

So far, we know that as long as $k\varepsilon^2 \to \infty$, we
have
$$
k^{-\alpha} L_k^\alpha(x/k) e^{-x/k} = (1+o(1))x^{-\alpha/2}
J_\alpha(2\sqrt{x}) + O\left(\sqrt{k}e^{-kC \varepsilon^2}\right)
+ O\left(\varepsilon k^{1/2-\delta/2}\right).
$$
Now if we take $\varepsilon = k^{-\beta}$ with $(1-\delta)/2 <
\beta < 1/2$, we find that
$$
k^{-\alpha} L_k^\alpha(x/k) e^{-x/k} = x^{-\alpha/2}
J_\alpha(2\sqrt{x}) + o(1),
$$
as desired.
\end{proof}

\begin{lemma}
\label{sumint}
For $\alpha > -1/2$, if $f \co [0,\infty) \to \R$ is
continuous and satisfies
$$
|f(x)| \le C(1+|x|)^{-\alpha-1-\varepsilon}
$$
for some $C>0$ and $\varepsilon>0$, then
$$
\sum_{k \ge 0} t^k \int_0^\infty f(x/y)L_k^\alpha(x) x^\alpha
e^{-x} \, dx = (1-t)^{-\alpha-1}\int_0^\infty f(x/y) x^\alpha
e^{-x/(1-t)} \, dx
$$
whenever $|t| < 1/3$.
\end{lemma}

\begin{proof}
We would like to convert this sum to
$$
\int_0^\infty \sum_{k \ge 0} f(x/y)L_k^\alpha(x) x^\alpha e^{-x}
t^k \, dx
$$
and apply the generating function
$$
\sum_{k \ge 0} L_k^\alpha(x) t^k = (1-t)^{-\alpha-1} e^{-xt/(1-t)}
$$
((6.2.4) from \cite{AAR}).  To do so, we must justify
interchanging the limit with the sum.

Let
$$
g(t) = (1-t)^{-\alpha-1} e^{-xt/(1-t)} = (1-t)^{-\alpha-1} e^x
e^{-x/(1-t)}.
$$
Then the Lagrange form of the remainder in Taylor's theorem implies
$$
g(t) = \sum_{k=0}^{m-1} L_k^\alpha(x) t^k + \frac{g^{(m)}(s)}{m!}
t^m
$$
for some $s$ satisfying $|s| \le |t|$.  By Lemma~\ref{induction},
$$
\frac{g^{(m)}(s)}{m!} = \pm e^x (1-s)^{-\alpha-1-m} e^{-x/(1-s)}
L_m^\alpha(x/(1-s)).
$$
It follows from Lemma~\ref{saddlept} that
$$
\left|e^{-x/(1-s)} L_m^\alpha(x/(1-s))\right| \le C' m^\alpha
$$
for some constant $C'>0$ (depending on $\alpha$).  Thus,
$$
\left|(1-t)^{-\alpha-1}\int_0^\infty f(x/y) x^\alpha e^{-x/(1-t)}
\, dx - \sum_{k = 0}^{m-1} t^k \int_0^\infty f(x/y)L_k^\alpha(x)
x^\alpha e^{-x} \, dx\right|
$$
is bounded above by
\begin{equation}
\label{upperbound}
C' \left(\int_0^\infty f(x/y) x^\alpha \, dx\right)
(1-s)^{-\alpha-1} m^\alpha \left(\frac{t}{1-s}\right)^m.
\end{equation}
The integral in \eqref{upperbound} is finite because of
the bound on $|f|$ in the lemma statement.
Because $|t| < 1/3$ and $|s| \le |t|$, we have
$$
\left|\frac{t}{1-s}\right| < \frac{1}{2},
$$
and hence \eqref{upperbound} tends to $0$ as $m \to \infty$.
\end{proof}

\begin{proposition}
\label{silpchar}
Let $\alpha > -1/2$, and suppose $f \co [0,\infty)
\to \R$ is continuous, and satisfies $|f(x)| \le
C(1+|x|)^{-\alpha-1-\varepsilon}$ for some $C>0$ and
$\varepsilon>0$. Then $f$ has the $\alpha$--SILP property iff
$$
f(x) = \int_0^\infty (xy)^{-\alpha/2} J_\alpha(2\sqrt{xy}) \,
dg(y)
$$
for some weakly increasing function $g$.
\end{proposition}

Note that one can compute directly the Laguerre coefficients of
the scalings of $x^{-\alpha/2} J_\alpha(2\sqrt{x})$ and verify
that they are positive (see Example~3 in Section~4.24
of~\cite{Leb}). Proposition~\ref{silpchar} tells us that this
function is essentially the only $\alpha$--SILP function.

\begin{proof}
We know that $f$ has the $\alpha$--SILP property iff for every
$y>0$,
$$
\sum_{k \ge 0} t^k \int_0^\infty f(x/y)L_k^\alpha(x) x^\alpha
e^{-x} \, dx
$$
has non-negative coefficients as a power series in $t$.  By
Lemma~\ref{sumint}, we can write this function (for small $t$)
as
$$
(1-t)^{-\alpha-1}\int_0^\infty f(x/y) x^\alpha e^{-x/(1-t)} \, dx,
$$
which is a positive constant (a power of $y$) times
$$
(1-t)^{-\alpha-1}\int_0^\infty f(x) x^\alpha e^{-xy/(1-t)} \, dx.
$$
Define $\tilde f$ to be the Laplace transform of $x \mapsto
x^\alpha f(x)$. Then $f$ has the $\alpha$--SILP property iff
$$
(1-t)^{-\alpha-1} {\tilde f}(y/(1-t))
$$
has non-negative coefficients as a power series in $t$. We can
rescale $t$ by a factor of $y$ and pull out a power of $y$ to see
that this happens iff
$$
(1/y-t)^{-\alpha-1} {\tilde f}(1/(1/y-t))
$$
has non-negative coefficients.  That happens for all $y>0$ iff the
function $u \mapsto u^{-\alpha-1} {\tilde f}(1/u)$ has successive
derivatives alternating in sign (the function is non-negative,
its derivative non-positive, its second derivative non-negative,
etc.). By Bernstein's theorem (Theorem~12b of Chapter~IV of
\cite[page 161]{W}), this holds iff it is the Laplace transform of
a positive measure.

Thus, we have shown that $f$ has the $\alpha$--SILP property iff
there is a weakly increasing function $g$ such that for $u>0$,
$$
u^{-\alpha-1} \int_0^\infty f(x) x^\alpha e^{-x/u} \, dx =
\int_0^\infty e^{-yu} \, dg(y).
$$

To finish proving the proposition, we can work as follows.  We
know that
$$
\int_0^\infty f(x) x^\alpha e^{-xu} \, dx =
u^{-\alpha-1}\int_0^\infty e^{-y/u} \, dg(y).
$$
We can now apply the following general theorem for inverting a
Laplace transform: if
$$
\phi(u) = \int_0^\infty \psi(x) e^{-xu} \, dx,
$$
then
$$
\psi(x) = \lim_{k \rightarrow \infty} \frac{(-1)^k}{k!}
\phi^{(k)}\left(\frac{k}{x}\right) \left(\frac{k}{x}\right)^{k+1}
$$
wherever $\psi$ is continuous.  (See Corollary~6a.2 of Chapter~VII
in \cite[page 289]{W}.)

We can apply this to our equation, and differentiate under the
integral sign (justified since the differentiated integrals
converge uniformly as $u$ ranges over any compact subset of
$(0,\infty)$; see Theorem~14 of Chapter~10 in \cite[page 358]{W2}).
Using Lemma~\ref{induction}, it follows that
$$
x^\alpha f(x) = \lim_{k \rightarrow \infty} \int_0^\infty
\left(\frac{k}{x}\right)^{-\alpha} L_k^\alpha
\left(\frac{xy}{k}\right) e^{-xy/k} \, dg(y).
$$
To finish the proof, we apply Lemma~\ref{saddlept}, but we need to
check that passage to the limit under the integral sign is
justified. Because of the uniform convergence, it is justified as
long as constant functions are integrable with respect to $dg$.
However, that is true, for the following reason.  By definition,
$g$ satisfies
$$
u^{-\alpha-1} \int_0^\infty f(x) x^\alpha e^{-x/u} \, dx =
\int_0^\infty e^{-yu} \, dg(y),
$$
which is equivalent to
$$
\int_0^\infty f(ux) x^\alpha e^{-x} \, dx =
\int_0^\infty e^{-yu} \, dg(y).
$$
When we let $u \rightarrow 0+$, the left side converges to
$$
f(0)\int_0^\infty x^\alpha e^{-x} \, dx
$$
(by the dominated convergence theorem: recall that $f$ is bounded
and continuous), so the right side converges as $u \rightarrow
0+$.  By monotone convergence, we see that constant functions are
integrable with respect to $dg$, which is what we need.
\end{proof}

\begin{corollary}
\label{silpcor} For integers $n\!>\!1$, a function $f\colon [0,\infty)
\!\to\! \R$ has the $(n/2-1)$--SILP property iff the function from
$\R^n$ to $\R$ given by $x \mapsto f(|x|^2)$ is continuous,
satisfies
$$
|f(|x|^2)| \le C(1+|x|)^{-n-\varepsilon}
$$
for some $C>0$ and $\varepsilon>0$, and is the Fourier transform
of a non-negative distribution.
\end{corollary}

Corollary~\ref{silpcor} follows from combining
Proposition~\ref{silpchar} with Theorem~9.10.3 of \cite{AAR}
(see Proposition~2.1
of \cite{CE}), after some changes of variables.  Using
Corollary~\ref{silpcor}, one can check with some simple
manipulations that for $n>1$, Theorem~\ref{silpbound} implies
Theorem~\ref{main} for lattice packings (and, as pointed out
above, the general case can be proved similarly). It is seemingly
more general, because it does not constrain the Fourier transform
at infinity. However, the additional generality does not seem
useful, and one could likely generalize the proof in \cite{CE} to
use a version of Poisson summation with fewer hypotheses (for example,
see Theorem~D.4.1 in \cite{AAR}).

\begin{corollary}
\label{productsilp}
For $\alpha > -1/2$, the product of two
$\alpha$--SILP functions is always an $\alpha$--SILP function.
\end{corollary}

Corollary~\ref{productsilp} follows immediately from
Corollary~\ref{silpcor} when $\alpha=n/2-1$ with $n \in \Z$, and
can be proved for arbitrary $\alpha$ using
Proposition~\ref{silpchar} together with 13.46~(3) of \cite{Wat}
or~(7) from~\S3 of \cite{Sz}. It seems surprisingly difficult to
prove directly from the definition of a SILP function:  it would
follow trivially if the product of two Laguerre polynomials were
a positive combination of Laguerre polynomials, but that is not
the case. In fact, the coefficients of such a product alternate
in sign; that is, the polynomials $(-1)^k L_k^\alpha$ have the
property that the set of positive combinations of them is closed
under multiplication.

\section{Optimality of Bessel functions}
\label{opt}

Let $j_\nu$ denote the first positive root of $J_\nu$.  According
to Proposition~6.1
of \cite{CE}, the function $f \co \R^n \to \R$ defined by
\begin{equation}
\label{levensh}
f(x) = \frac{J_{n/2}(j_{n/2}|x|)^2}{(1-|x|^2)|x|^n}
\end{equation}
satisfies the hypotheses of Theorem~\ref{main}, and leads to
the upper bound
$$
\frac{j_{n/2}^n}{(n/2)!^2 4^n}
$$
for the densities of $n$--dimensional sphere packings.  The
Fourier transform $\widehat{f}$ has support in the ball of radius
$j_{n/2}/\pi$ about the origin.  We will show that among all such
functions, $f$ proves the best sphere packing bound.  This is
analogous to a theorem of Sidel'nikov \cite{S} for the case of
error-correcting codes and spherical codes.  It was first proved
in the setting of sphere packings by Gorbachev \cite{G}.  Our
proof will be based on the same identity as Gorbachev's, but the
proof of the identity appears to be new.

For notational simplicity, we view $f$ and $\widehat{f}$ as
functions on $[0,\infty)$; that is, $f(r)$ will denote the common
value of $f$ on all vectors of length $r$. Let $\nu = n/2-1$, and
let $\lambda_1 < \lambda_2 < \cdots$ be the positive roots of
$J_{\nu+1}(x)$ (equivalently, the positive roots of $-\nu
J_\nu(x) + xJ_\nu'(x)$;  see equation (4) in \S 3.2 of
\cite{Wat}). Define $B_r(x)$ to be the closed ball of radius $r$
about $x$.

Our main technical tool is the following identity due to Ben
Ghanem and Frappier (the $p=0$ case of Lemma~4 in \cite{Be}), who
state it with weaker technical hypotheses and a different proof.

\begin{theorem}[Ben Ghanem and Frappier \cite{Be}]
\label{dualthm} Let $f \co \R^n \to \R$ be a radial Schwartz
function.  If $\supp(\widehat{f}\,) \subseteq B_r(0)$, then
$$
\widehat{f}(0) = \frac{(n/2)!2^{n}}{\pi^{n/2}r^n}f(0) +
\sum_{m=1}^\infty \frac{4 \lambda_m^{n-2}}{(n/2-1)! \pi^{n/2} r^n
J_{n/2-1}(\lambda_m)^2} f\left(\frac{\lambda_m}{\pi r}\right).
$$
\end{theorem}

We will postpone the proof of Theorem~\ref{dualthm} until we have
developed several lemmas.  First, however, we deduce the desired
optimality:

\begin{corollary}[Gorbachev \cite{G}]
\label{optcor} Suppose $f \co \R^n \rightarrow \R$ is a radial,
admissible function, is not identically zero, and satisfies the
following three conditions:
\begin{enumerate}
\item $f(x)\le 0$ for $|x| \ge 1$,

\item ${\widehat f}(t)\ge 0$ for all $t$, and

\item $\supp(\widehat{f}\,) \subseteq B_{j_{n/2}/\pi}(0)$.
\end{enumerate}
Then
$$
\frac{\pi^{n/2}}{(n/2)!2^n} \cdot \frac{f(0)}{\widehat{f}(0)} \ge
\frac{j_{n/2}^n}{(n/2)!^2 4^n}.
$$
\end{corollary}

\begin{proof}[Proof of Corollary~\ref{optcor}]
Let $r = j_{n/2}/\pi$. If $f$ were a Schwartz function, then
Theorem~\ref{dualthm} would imply that
$$
\widehat{f}(0) \le
\frac{(n/2)!2^{n}}{\pi^{n/2}(j_{n/2}/\pi)^n}f(0),
$$
since $\lambda_m/(\pi r) \ge 1$ for $m\ge1$.  For more general
functions $f$, the series
$$
\frac{(n/2)!2^{n}}{\pi^{n/2}r^n}f(0) + \sum_{m=1}^\infty \frac{4
\lambda_m^{n-2}}{(n/2-1)! \pi^{n/2} r^n J_{n/2-1}(\lambda_m)^2}
f\left(\frac{\lambda_m}{\pi r}\right)
$$
at least still converges, since the terms are
$O(m^{-1-\varepsilon})$ for some $\varepsilon>0$ (namely, the
$\varepsilon$ from the definition of admissibility); to verify
this, note that $\lambda_m$ grows linearly with $m$, and that
$J_\nu(z)^2 + J_{\nu+1}(z)^2 \sim 2/(\pi z)$ (see Section~7.21 of
\cite[page 200]{Wat}), so $J_{\nu}(\lambda_m)^2 \sim 2/(\pi
\lambda_m)$.  However, we must verify that it converges to
$\widehat{f}(0)$.

We need to smooth $\widehat{f}$ without increasing its support.
Let $i_\delta$ denote any non-negative, smooth function of
integral~$1$ with support in the ball of radius $\delta$ about
the origin. Let $f_\varepsilon(x) =
f(x(1-\varepsilon))\widehat{\imath}_{r\varepsilon/2}(x)$, where
$r = j_{n/2}/\pi$.  This is a Schwartz function whose Fourier
transform has support in the ball of radius $r(1-\varepsilon/2)$,
so Theorem~\ref{dualthm} applies to $f_\varepsilon$. As
$\varepsilon \to 0+$, the functions $f_\varepsilon$ and
$\widehat{f}_\varepsilon$ converge pointwise to $f$ and
$\widehat{f}$, respectively.  Since
$|\,\widehat{\imath}_{r\varepsilon/2}| \le 1$ everywhere,
dominated convergence lets us interchange the limit as
$\varepsilon \to 0+$ with the sum over $m$ to conclude that
$$
\widehat{f}(0) = \frac{(n/2)!2^{n}}{\pi^{n/2}r^n}f(0) +
\sum_{m=1}^\infty \frac{4 \lambda_m^{n-2}}{(n/2-1)! \pi^{n/2} r^n
J_{n/2-1}(\lambda_m)^2} f\left(\frac{\lambda_m}{\pi r}\right),
$$
and we finish the proof as before.
\end{proof}

\begin{lemma}
\label{interp} Let $f \co \R^n \to \R$ be a radial Schwartz
function.  If $\supp(\widehat{f}\,) \subseteq B_r(0)$, then
for $u \in [0,1)$,
$$
2\pi \widehat{f}(ru) r^{\nu+2} = \frac{2 \Gamma(\nu+2)}{(\pi
r)^\nu} f(0) + \sum_{m=1}^\infty \frac{2 (\lambda_m/(2\pi r))^\nu
f(\lambda_m/(2\pi r))} {J_\nu(\lambda_m)^2} \frac{J_\nu(\lambda_m
u)}{u^\nu}.
$$
The same holds even if $\widehat{f}$ is not smooth at radius $r$
(but is left continuous at radius $r$, and still smooth at all
smaller radii), as long as the values of $f$ in the sum decrease
faster than any power of $1/m$ as $m \to \infty$.
\end{lemma}

Note that if $f$ is a Schwartz function, then the condition on
the decay of the values of $f$ automatically holds.

\begin{proof}
Because $\supp(\widehat{f}\,) \subseteq B_r(0)$,
we have
$$
x^\nu f(x) = \int_0^1 g(u) u^{\nu+1} J_\nu(2\pi r u x) \, du,
$$
where $g(u) = 2\pi \widehat{f}(ru) r^{\nu+2}$ (see
Theorem~9.10.3 of \cite{AAR}, or
Proposition~2.1
of \cite{CE}). We begin by expanding $g(u)u^\nu$ into a Dini
series. For a quick introduction to Dini series, see
\cite[page 130]{Leb}.  Unfortunately, for a technical reason that
reference does not cover the case we need here (see footnote~33 on
page~130).  For a more thorough reference, which covers
everything we need, see Sections~18.3--18.35 of \cite{Wat}. In
Watson's notation, we are dealing with the case $H+\nu=0$ (see
page 597 of \cite{Wat}).  Convergence of the Dini series to
$g(u)u^\nu$ for $u \in (0,1)$ follows from standard results (see
pages 601--602 of \cite{Wat}), and at $u=0$ it follows from
continuity of $g$ at $0$ and uniform convergence of the Dini
series (which itself follows from the decay of $f(\lambda_m/(\pi
r))$).

The Dini series expansion of $g(u) u^{\nu}$ is
$$
g(u) u^{\nu} = 2(\nu+1)u^\nu \int_0^1 t^{\nu+1} g(t) t^\nu \, dt
+ \sum_{m=1}^\infty b_m J_\nu(\lambda_m u),
$$
where
\begin{eqnarray*}
b_m &=& \frac{2\lambda_m^2}{(\lambda_m^2-\nu^2)J_\nu(\lambda_m)^2
+ \lambda_m^2 J_\nu'(\lambda_m)^2}
\int_0^1 t g(t) t^\nu J_\nu(\lambda_m t) \, dt\\
&=& \frac{2\lambda_m^2 (\lambda_m/(2\pi r))^\nu}
{(\lambda_m^2-\nu^2)J_\nu(\lambda_m)^2 + \lambda_m^2
J_\nu'(\lambda_m)^2} f(\lambda_m/(2\pi r)).
\end{eqnarray*}
Note also that
$$
\lim_{x \to 0} \int_0^1 g(u) u^{\nu+1} J_\nu(2\pi r u x)/x^\nu \,
du = \int_0^1 g(u) u^{\nu+1} \frac{(\pi r
u)^\nu}{\Gamma(\nu+1)}\, du,
$$
since as $x \to 0$,
$$
\frac{J_\nu(x)}{x^\nu} \to \frac{1}{2^\nu \Gamma(\nu+1)},
$$
so
$$
\int_0^1 t^{\nu+1} g(t) t^\nu \, dt = f(0) \Gamma(\nu+1) / (\pi
r)^\nu.
$$
Furthermore, $\lambda_m J_\nu'(\lambda_m) = \nu
J_\nu(\lambda_m)$, so
$$
(\lambda_m^2-\nu^2)J_\nu(\lambda_m)^2 + \lambda_m^2
J_\nu'(\lambda_m)^2 = \lambda_m^2 J_\nu(\lambda_m)^2.
$$
Thus,
$$
g(u) = \frac{2 \Gamma(\nu+2)}{(\pi r)^\nu} f(0) +
\sum_{m=1}^\infty \frac{2 (\lambda_m/(2\pi r))^\nu
f(\lambda_m/(2\pi r))} {J_\nu(\lambda_m)^2} \frac{J_\nu(\lambda_m
u)}{u^\nu},
$$
as desired.
\end{proof}

\begin{lemma}
\label{pws} Let $f$ be a function from $[0,\infty)$ to $\R$.  The
function $x \mapsto f(|x|)$ from $\R^n$ to $\R$ is the Fourier
transform of a compactly support distribution iff $f$ extends to
an even, entire function on $\C$ that satisfies
$$
|f(z)| \le C (1+|z|)^k e^{C' | \Imag z|}
$$
for some $C$, $C'$, and $k$.
\end{lemma}

\begin{proof}
This lemma is essentially a special case of the
Paley-Wiener-Schwartz theorem (Theorem~7.3.1 in \cite{H}).  The
only difference is that the general theorem is not restricted to
radial functions, and characterizes Fourier transforms of
compactly supported distributions as entire functions $g$ of $n$
complex variables satisfying
\begin{equation}
\label{pwsbound}
|g(z_1,\dots,z_n)| \le
C\left(1+\sqrt{|z_1|^2+\dots+|z_n|^2}\right)^k e^{C'\sqrt{(\Imag
z_1)^2 + \dots +(\Imag z_n)^2}}.
\end{equation}
The only subtlety in deriving the lemma from the general theorem
is in showing that if $f$ satisfies the hypotheses above, then
the function $g$ defined by
$$
g(z_1,\dots,z_n) = f\left(\sqrt{z_1^2 + \dots + z_n^2}\right)
$$
satisfies \eqref{pwsbound}.  To do that, the elementary inequality
$$
\left|\Imag \sqrt{z_1^2 + \dots + z_n^2}\right| \le \sqrt{(\Imag
z_1)^2 + \dots +(\Imag z_n)^2}
$$
can be used.  To prove that inequality, one can use induction to
reduce to the $n=2$ case, and prove that case by direct
manipulation of both sides.
\end{proof}

Now we are ready to prove Theorem~\ref{dualthm}.  Notice that it
says that to determine the integral of $f$, we need only half as
many values as we need to reconstruct the whole function via
Lemma~\ref{interp}.  This phenomenon is analogous to Gauss-Jacobi
quadrature (see Theorem~14.2.1 of \cite{D}). The proof given
below is in fact modeled after the proof of Gauss-Jacobi
quadrature, although carrying it out rigorously is more involved.

\begin{proof}[Proof of Theorem~\ref{dualthm}]
Let $\varepsilon>0$, and
define $\tilde h \co [-1,1] \to \R$ by
$$
\tilde h(u)  = \frac{2 \Gamma(\nu+2)}{(\pi
(r/2+\varepsilon))^\nu} f(0) + \sum_{m=1}^\infty \frac{2
\left(\frac{\lambda_m}{2\pi (r/2+\varepsilon)}\right)^\nu
f\left(\frac{\lambda_m}{2\pi (r/2+\varepsilon)}\right)}
{J_\nu(\lambda_m)^2} \frac{J_\nu(\lambda_m u)}{u^\nu}.
$$
(The functions $J_\nu(\lambda_m u)/u^\nu$ are even, so this is no
different from defining $\tilde h$ on $[0,1]$.) Since $f$ is a
Schwartz function, the values of $f$ in the series above decrease
quickly enough that it defines a $C^\infty$ function on $(-1,1)$.
Define $\widehat{h}$ by
$$
2\pi \widehat{h}((r/2+\varepsilon)u) (r/2+\varepsilon)^{\nu+2} =
\begin{cases}
\tilde h(u) & \textup{if $|u| \le 1$, and}\\
0 & \textup{otherwise,}
\end{cases}
$$
and define $h$ to be the Fourier transform of $\widehat{h}$. Then
$\supp(\widehat{h}) \subseteq B_{r/2+\varepsilon}(0)$. By
Lemma~\ref{interp}, combined with uniqueness for Dini series
(which follows from orthogonality), we have
$$
h\left(\frac{\lambda_m}{2\pi (r/2+\varepsilon)}\right)
=
f\left(\frac{\lambda_m}{2\pi (r/2+\varepsilon)}\right)
$$
for all $m$, and $h(0)=f(0)$.  (Note that $\widehat{h}$ may not
be smooth at radius $r/2+\varepsilon$, but that does not violate
the hypotheses of Lemma~\ref{interp}.)

Now let $\newchi_R$ denote the characteristic function of a ball
of radius~$R$ about the origin, so that
$$
\newchihat_R(x) = J_{n/2}(2\pi R |x|) (R/|x|)^{n/2}.
$$
The entire function $f-h$ has roots wherever
$\newchihat_{r/2+\varepsilon}$ does, and
$\newchihat_{r/2+\varepsilon}$ has only single roots, so the
quotient $g = (f-h)/\newchihat_{r/2+\varepsilon}$ is entire.

We would like to conclude that $g$ is the Fourier transform of a
compactly supported distribution.  By Lemma~\ref{pws}, this
requires bounds for $g$, and it is not obvious that dividing by a
Bessel function does not ruin the bounds.  We prove this in two
steps.  First, Lemma~1 of \cite{Lev} implies (after rescaling
variables) that
$$
|J_{n/2}(z)/z^{n/2}| \ge \frac{c_1 e^{c_2 |\Imag
z|}}{(1+|z|)^{c_3}}
$$
whenever $|\Imag z| \ge c_4$, for some constants
$c_1,c_2,c_3,c_4$, with $c_1>0$ of course. That means that
dividing by it does not mess up our bounds when the absolute
value of the imaginary part is at least $c_4$.  The second step
is to deal with points near the real axis.  Consider a box with
sides on the lines with imaginary part $\pm c_4$ and real part
$\pm(k\pi + (\pi n + 1)/4)$, where $k$ is a positive integer.  By
the maximum principle, the maximum of $g$ over the interior of the
box must occur on the sides.  We know that $g$ satisfies the
bound we want on the top and bottom, and $g$ is even, so we only
need to estimate $g$ on the right side.

For $z$ in the right half-plane, we have
\begin{eqnarray*}
J_{n/2}(z) = \sqrt{\frac{2}{\pi z}} \left(\cos\left(z -
\frac{\pi n +1}{4}\right)(1+O(1/z^2))\right.\\
\left.\phantom{} + \sin\left(z - \frac{\pi n
+1}{4}\right)(O(1/z))\right)
\end{eqnarray*}
(see~(1) in Section~7.21 of \cite{Wat}).  When $z$ has real part
$k\pi$, we have $\cos(z) = (-1)^k \cosh(\Imag z)$, which has
absolute value at least $1$.  Thus, on the right side of the box,
the cosine factor is always at least $1$.  The sine factor is
bounded, because $\Imag z$ is bounded, so we see that on the
right side of the box $J_{n/2}(z)/z^{n/2}$ is never smaller than
a power of $1/|z|$.

When we combine these estimates, it follows from Lemma~\ref{pws}
that $g$ is the Fourier transform of a distribution with compact
support. Furthermore, the Titchmarsh-Lions theorem (see
Theorem~4.3.3 in \cite{H}) implies that the convex hull of the
support of $f-h$ equals the Minkowski sum of those of
$\widehat{g}$ and $\newchi_{r/2+\varepsilon}$, so
$\supp(\widehat{g}) \subseteq B_{r/2-\varepsilon}(0)$.

Let $i_\delta$ denote any non-negative, smooth function of
integral~$1$ with support in the ball of radius $\delta$ about
the origin. We have
$$
f \widehat{\imath}_\delta - h \widehat{\imath}_\delta =
(\newchihat_{r/2+\varepsilon} \widehat{\imath}_\delta) g.
$$
Now both sides are integrable functions (note that this is not
obviously true of either $h$ or $g\newchihat_{r/2+\varepsilon}$,
which is why we had to multiply by $\widehat{\imath}_\delta$),
and we find that
$$
(\widehat{f} * i_\delta)(0) - (\widehat{h} * i_\delta)(0) =
\int (\newchi_{r/2+\varepsilon} * i_\delta)\widehat{g} .
$$
Because $\supp(\widehat{g}) \subseteq B_{r/2-\varepsilon}(0)$, if
we take $\delta < \varepsilon$ we have
$$
\int  (\newchi_{r/2+\varepsilon} * i_\delta) \widehat{g} =
\int \widehat{g} = g(0) = 0,
$$
where $g(0)=0$ because $f(0)=h(0)$.  Thus.
$$
(\widehat{f} * i_\delta)(0) = (\widehat{h} * i_\delta)(0).
$$
If we let $\delta \to 0+$, we find that
$\widehat{f}(0)=\widehat{h}(0)$, because both $\widehat{f}$ and
$\widehat{h}$ are continuous near~$0$. It follows from the way
$\widehat{h}$ was defined that $\widehat{f}(0)$ equals
$$
\frac{(n/2)!2^{n}}{\pi^{n/2}(r+2\varepsilon)^n}f(0) +
\sum_{m=1}^\infty \frac{4 \lambda_m^{n-2}}{(n/2-1)! \pi^{n/2}
(r+2\varepsilon)^n J_{n/2-1}(\lambda_m)^2}
f\left(\frac{\lambda_m}{\pi (r+2\varepsilon)}\right).
$$
Now sending $\varepsilon \to 0+$ proves the desired result, by
dominated convergence.
\end{proof}

\section{The dual program}
\label{dualprog}

It is natural to view choosing the optimal function $f$ in
Theorem~\ref{main} as solving an infinite-dimensional linear
programming problem:  if we fix $\widehat{f}(0)=1$, then we are
trying to minimize the linear functional $f(0)$ of $f$, subject to
linear inequalities on $f$.  The technicalities are slightly
subtle;  for example, it is not immediately clear what the right
space of functions to consider is (admissibility might be too ad
hoc). It seems likely that Schwartz functions suffice.  One can
come arbitrarily close to the optimum with functions $f$ such that
$f$ and $\widehat{f}$ are smooth and rapidly decreasing, where we
say $g \co \R^n \to \R$ is rapidly decreasing if $g(x) =
O((1+|x|)^{-k})$ for every $k>0$: given any $f$ that satisfies
the hypotheses of Theorem~\ref{main}, let
$$
f_\varepsilon(x) = ((f * i_\varepsilon * i_\varepsilon)
{\widehat{\imath}}_\varepsilon\,\!\!^2)((1+2\varepsilon)x),
$$
where $i_\varepsilon$ is any non-negative, smooth function of
integral~$1$ with support in $B_\varepsilon(0)$.  Then
$f_\varepsilon$ has the desired properties, still obeys the
required inequalities, and satisfies
$$
\lim_{\varepsilon \to
0}\frac{f_\varepsilon(0)}{\widehat{f_\varepsilon}(0)} =
\frac{f(0)}{\widehat{f}(0)}.
$$
Presumably Schwartz functions also come arbitrarily close, but
one would have to worry about making the derivatives rapidly
decreasing as well.  Despite the fact that rapidly decreasing
functions come close to the optimal bounds, it is not clear
whether they reach them. For example, even for $n=1$, where one
can write down several explicit functions that solve the sphere
packing problem (see Sections~3 and~5 of \cite{CE}),
these functions are not rapidly decreasing.

In this context, it is natural to study the dual linear program,
to prove bounds on how good the sphere packing bounds produced by
Theorem~\ref{main} can be. The results of Section~\ref{opt}
amount to doing exactly this, for a restricted linear program in
which we limit the support of $\widehat{f}$. Unfortunately, in
the unrestricted case the dual program seems no easier to solve
in general than the primal program is. However, it leads to
several intriguing open problems.

One formulation of the dual program is as follows: find the
largest $c$ such that there is a tempered distribution $g$ on
$\R^n$ satisfying
\begin{enumerate}
\item $g = \delta + h$ with $h \ge 0$,

\item $\supp(h) \subseteq \{x : |x|\ge 1\}$, and

\item $\widehat{g} \ge c \delta$.
\end{enumerate}
Here $\delta$ is a delta function at the origin, and inequalities
between distributions mean that applying both sides to
non-negative functions preserves this inequality.  For $g$
satisfying~(1)--(3) above, and any radial function $f$ satisfying
the hypotheses of Theorem~\ref{main} such that $f$ and
$\widehat{f}$ are rapidly decreasing, we have
$$
f(0) \ge \int_{\R^n} fg = \int_{\R^n} \widehat{f}\,\widehat{g}
\ge c \widehat{f}(0).
$$
Here, we use the fact that one can apply a non-negative tempered
distribution to any rapidly decreasing function, because
non-negative tempered distributions are exactly measures $\mu$
such that
$$
\int_{\R^n} \frac{d\mu(x)}{(1+|x|)^k} < \infty
$$
for some $k$ (see Theorem VII in Chapter 7, \S4 of
\cite[page 242]{Sch}). Thus $f(0)/\widehat{f}(0) \ge c$. The duality
theorem of linear programming suggests that there is no gap
between the smallest $f(0)/\widehat{f}(0)$ and largest $c$, but
it is not clear how to prove it in this infinite-dimensional
setting.

Given any lattice $\Lambda$ with minimum non-zero vector
length~$1$, summing over $\Lambda$ defines a tempered
distribution that clearly satisfies properties~(1) and~(2), and
Poisson summation implies that it has property~(3) as well.  As
is the case for the functions $f$, we can rotationally symmetrize
$g$, so that  $g$ and $\widehat{g}$ are positive linear
combinations of spherical delta functions, where we define a
spherical delta function $\delta_r$ on $\R^n$ to be a
distribution with support on the sphere of radius $r$ about the
origin, such that integrating any function times $\delta_r$ gives
the average of that function over the sphere. One would expect
that the optimal radial $g$ should always be a linear combination
of spherical delta functions, but it is not clear how to prove it.
Aside from the origin, $g$ and $\widehat{g}$ should be supported
on the zeros of the optimal $f$ and $\widehat{f}$, respectively,
but why must these zeros even occur at a discrete set of radii?

\begin{question}
\label{lindepq}
Consider tempered distributions $g$ such that $g$
and $\widehat{g}$ are linear combinations of spherical delta
functions. Is every such distribution in the span of the
rotationally symmetrized Poisson summation distributions?
\end{question}

It seems very unlikely that the answer to Question~\ref{lindepq}
is yes. Any counterexample would be of interest, since the
optimal distributions $g$ in most dimensions (not $1$, $2$, $8$,
or $24$) are probably counterexamples.

One interesting case is $72$ dimensions.  It is an open question
whether there exists an ``extremal lattice of Type~II'' in $\R^{72}$,
in other words, an even unimodular lattice in $\R^{72}$ with minimal
non-zero norm at least $8$ (see \cite[page 194]{SPLAG} for more
details).  Such a lattice might be as extraordinary as $E_8$ or the
Leech lattice. Unfortunately, it seems unlikely that one exists.
However, its existence cannot be ruled out by Theorem~\ref{main}. The
simplest way to see that is in light of Section~\ref{newproof}. A
proof that the lattice did not exist would amount to a proof that its
theta series could not exist.  However, although the extremal lattice
may not exist, there is a modular form that would be its theta series
if it did exist (see \cite[page 195]{SPLAG}).  In fact, the modular
form comes from a distribution $g$ as above, because it is a
polynomial in the theta series of $E_8$ and the Leech lattice, and
therefore comes {}from a $g$ that is the corresponding linear
combination of Poisson summation for direct sums of $E_8$ and the
Leech lattice.  If $\Theta_n$ denotes the theta series of $E_8$, the
Leech lattice, and the hypothetical $72$--dimensional lattice for $n
=8,24,72$, respectively, then
$$
\Theta_{72} = \frac{79}{1080} \Theta_{24}^3 +
\frac{1183}{720}\Theta_{24}^2\Theta_8^3 -
\frac{91}{180}\Theta_{24}\Theta_8^6 - \frac{91}{432}\Theta_8^9.
$$
Despite the minus signs, all the coefficients of $\Theta_{72}$
are non-negative.

The most elegant form of the dual program comes from a rescaling
analogous to that in Theorem~3.2
of \cite{CE}.  Define a \textit{relaxed lattice\/}
to be a tempered
distribution $g$ such that $g$ and $\widehat{g}$ are of the form
$$
\sum_{i \ge 0} a_i \delta_{r_i}
$$
with $a_i \ge 0$ for all $i$ (not all $0$), and $0 = r_0 < r_1 <
r_2 < \cdots$.  Call a relaxed lattice $g$ \textit{self-dual\/}
if $\widehat{g}=g$.  How large can $r_1$ be?

\begin{conjecture}
\label{dualconj}
In every dimension, the largest possible value
of $r_1$ in a self-dual relaxed lattice equals the smallest value
of $r$ possible in Theorem~3.2
of \cite{CE}.
\end{conjecture}

One might imagine that the self-duality in Conjecture~\ref{dualconj}
would follow from some sort of symmetry of the linear programming
problem, but that is not clear.  If this conjecture is true, it would
explain the otherwise remarkable fact that the minimal values of $r$
in Proposition~7.1 and Theorem~3.2 of \cite{CE} always seem to agree
(see Conjecture~7.2 in that paper).

\end{document}